\documentclass[12pt]{article}
\usepackage{amsfonts}
\usepackage[english]{babel}
\usepackage[T1]{fontenc}
\usepackage{amsthm,amsmath,amsfonts,amssymb,epsfig,latexsym}
\usepackage{bm}
\usepackage{helvet}         
\usepackage{courier}        
\usepackage{multicol}        

\makeindex             

\newcommand{\ve}{\varepsilon}

\newcommand{\NN}{\mathbb{N}}

\newcommand{\ff}{\hat f}

\newcommand{\R}{\mathbb R}
\def\be#1\ee{\begin{equation}#1\end{equation}}
\newcommand{\fer}[1]{(\ref{#1})}

\newcommand{\bq}{\begin{equation}}
\newcommand{\eq}{\end{equation}}

\def\bqa{\begin{eqnarray}}
\def\eqa{\end{eqnarray}}

\def\e{\epsilon}

\def\ee{\varepsilon}

\newcommand{\bd}{\begin{displaymath}}
\newcommand{\ed}{\end{displaymath}}
\newcommand{\ba}{\begin{eqnarray}}
\newcommand{\ea}{\end{eqnarray}}

\def\ff{\widehat f}

\def\N{\mathbb{N}}
\def\R {{I \!\! R}}

\newenvironment{equations}{\equation\aligned}{\endaligned\endequation}

\begin{document}

\title{A Rosenau-type approach to the approximation \\ of the  linear Fokker--Planck equation}

\author{G. Toscani\thanks{Department of Mathematics, University of Pavia, and and IMATI-CNR,  via Ferrata 1, 27100 Pavia, Italy.
\texttt{giuseppe.toscani@unipv.it}} }
\date{}

\maketitle
\noindent
{\bf Abstract:} \small{The numerical approximation of the solution of the Fokker--Planck equation is a challenging problem that has been extensively investigated starting from the pioneering paper of Chang and Cooper  in 1970 \cite{CC70}. We revisit this problem at the light of the  approximation of the solution to the heat equation proposed by Rosenau \cite{Ro92}. Further, by means of the same idea, we address the problem of a consistent approximation to higher-order linear diffusion equations.}
\vskip 5mm

\noindent
{\bf Keywords}:
{Fokker--Planck equation; discrete schemes; Wild sums; Fourier-based me\-tri\-cs; higher-order diffusions.}


\section{Introduction}
The Fokker-Planck equation is a partial differential equation describing the time evolution of a density function $f(v,t)$, where $v\in \R^d, \, d \ge 1$ and $t \ge 0$, departing from a nonnegative initial density $\varphi(v)$. The standard assumptions on $\varphi(v)$ is that it possesses finite mass $\rho$, mean velocity $u$ and temperature $\theta$, where for any given density $g(v)$
 \be\label{mass}
 \rho = \int_{\R^d} g(v) \, dv
 \ee
is the mass density,
\be\label{velo}
 u = \frac 1\rho \int_{\R^d}v g(v) \, dv
 \ee
is the mean velocity, and $\theta$ is the temperature defined by
 \be\label{temp}
 \theta= \frac 1{d\rho} \int_{\R^d}|v-u|^2 g(v) \, dv.
 \ee
 The Fokker--Planck equation is a fundamental model in  kinetic theories and statistical mechanics. Its general form reads
 \be\label{FPgen}
\frac{\partial f}{\partial t} = J_{FP}(f) = \gamma \sum_{k=1}^d \left\{\frac{\partial^2 f}{\partial v_k^2} + \frac 1\theta  \frac{\partial }{\partial v_k}[(v_k-u_k)f]\right\}.
\ee
The one-particle friction constant $\gamma$ is usually assumed to be a function of $\rho, u, \theta$. Equation \fer{FPgen} has a stationary solution of given mass $\rho$, mean velocity $u$ and temperature $\theta$ given by the Maxwellian density function
 \be\label{maxw}
\mathcal{M}_{\rho,u,\theta}(v) = \rho \, \frac 1{(2\pi\theta)^{d/2}}\exp\left\{-\frac{|v-u|^2}{2\theta}\right\},
 \ee
which is such that $J_{FP}(M_{\rho,u,\theta})=0$.  The Fokker-Planck equation appears in many different contexts. It was originally derived for the distribution function of a Brownian particle in a fluid \cite{Ch43}, and is applicable in a more general form to a plasma \cite{ChC58}. A detailed investigation of this model has been performed by Frisch, Helfand, and Lebowitz \cite{FHL60} in connection with the kinetic theory of liquids. As shown more recently \cite{Tos98} (cf. also \cite{Cer88}), it provides also a good description of the grazing collisions in a one-dimensional gas. The Fokker-Planck operator $J_{FP}$  has the usual conservation properties of mass, mean velocity, and temperature, and $\int log f J_{FP}(f)\, dv <0$, which guarantees the increasing in time of the (Shannon) entropy
 \be\label{Sha}
 H(f)(t) = - \int_{\R^d} f(v,t) \log f(v,t) \, dv.
 \ee
It is interesting to remark that, if the friction $\gamma$ is taken to be proportional to the pressure $p = \rho\theta$ , then $J_{FP}(f)$ has the same kind of nonlinearity (quadratic) as the true Boltzmann equation.
 
For the purpose of accurate numerical simulations, a discretized Fokker-Planck equation must guarantee most of the conservation laws of the original equation, starting from mass conservation.  Furthermore, since the solution of the Fokker--Planck equation represents a  density function, any numerical scheme that approximates  equation \fer{FPgen} is required to guarantee the positivity of the solution. In addition, it would be desirable that an approximation scheme must be accurate and stable. 

The seminal paper for the approximation to equation \fer{FPgen} is due to Chang and Cooper \cite{CC70}. Other classical references are the paper by Larsen, Levermore Pomraning and Sanderson \cite{LLPS85}, and the well-known book by Risken \cite{Ri89}. 
Various aspects of the numerical approximation of Fokker--Planck equation were subsequently dealt with by a number of authors \cite{BD,BDS, De1, De2, Epp, MK}. Also in recent times, this problem has attracted the interest of research \cite{MB, PZ}. 

The aim of this paper is to present a discretized version of equation \fer{FPgen} which maintains most of the physical properties of the original equation. These properties include conservation of mass and positivity of the discrete solution, same evolution for the mean velocity and temperature, monotonicity in time of the discrete Shannon entropy, and the existence of an explicit discrete equilibrium density. In addition, the problem of the large-time behavior of the approximation and the convergence to the corresponding equilibrium density has been dealt with in the one-dimensional situation.

 This discrete version is largely inspired by a recent paper \cite{RT}, in which the kinetic meaning of the approximation to the heat equation proposed by Rosenau in \cite{Ro92} has been deeply investigated. 

The last part of the paper is devoted to show how this idea could be fruitfully applied to construct a numerical approximation to one-dimensional linear diffusion equations of higher order. In particular, it is shown that starting from this approximation one can easily obtain an explicitly computable formula for the central difference approximation of a derivative of any even order. 

\section{Main properties of the Fokker--Planck equation}\label{model}
\setcounter{equation}{0}

Given a nonnegative initial value $\varphi(v)$  with finite mean velocity $u$ and temperature $\theta$, easy computations show that the mass, mean velocity, and temperature of the solution to  the Fokker--Planck equation \fer{FPgen} do not change with time. It is convenient to normalize $f$ to be a probability density instead of a mass density, and change equation \fer{FPgen} to a dimensionless form. To do this, one introduces the dimensionless variables $\bar v$, $\bar t$, and the dimensionless functions $\bar\varphi$, $\bar f$ defined by the formulas
 \begin{equations}\label{dim}
& \bar v = \frac{v-u}{\sqrt\theta}, \qquad \bar t = \frac\gamma\theta t, \\
& \varphi(v) = \rho \theta^{-d/2}\bar\varphi(\bar v), \quad f(v,t) =\rho \theta^{-d/2}\bar f(\bar v, \bar t). 
\end{equations}
Substituting \fer{dim} into \fer{FPgen}, carrying out elementary calculations, and then omitting the bars, we obtain that the function $f(v, t)$ will now satisfy the equation
 \be\label{FP}
\frac{\partial f}{\partial t} = \bar J_{FP}(f) = \sum_{k=1}^d \left\{\frac{\partial^2 f}{\partial v_k^2} +  \frac{\partial }{\partial v_k}(v_k f]\right\}.
\ee
with the initial condition $\varphi(v)$ and consequently $f(v, t )$ satisfying the following simple normalization conditions
 \be\label{norm}
 \rho=1, \qquad u=0, \qquad \theta=1.
 \ee
The normalization \fer{norm} corresponds to the equilibrium Maxwellian density
\be\label{max}
 \mathcal{M}(v) =  \frac 1{(2\pi)^{d/2}}\exp\left\{-\frac{|v|^2}{2}\right\}.
 \ee
Let $\varphi$  be any probability density on $\R^d$ with finite second moment. Let $X$ be any random variable with this density, and let $W$ be any independent Gaussian random variable with density $M$ given by \fer{max}. For every $t > 0$ define
 \be\label{sol}
 Z_t = e^{-t}X + (1-e^{-2t})^{1/2}W.
 \ee
Then, the random variable $Z_t$ has a density $f(v,t)$ at each $t \ge 0$, and it is well-known that $f(t)$ is evolved from $\varphi$ under the action of the adjoint Ornstein--Uhlenbeck semigroup. Therefore $f(v,t)$ satisfies equation \fer{FP}, which can of course be checked directly from the definition.

Mean velocity and temperature of the solution at any time $t \ge 0$ can be obtained directly from expression \fer{sol}. We obtain
 \be\label{me1}
 u(t) = \langle Z_t\rangle = e^{-t}\langle X \rangle = e^{-t}\int_{\R^d} v\, \varphi(v) \, dv,
 \ee 
and
 \begin{equations}\label{te1}
 \theta(t) = \langle |Z_t|^2\rangle =&\,\, e^{-2t}\langle |X|^2 \rangle +  \left(1-e^{-2t}\right)\langle |M|^2 \rangle=\\
& 1- e^{-2t}\left( 1- \int_{\R^d} |v|^2\, \varphi(v) \, dv\right).
 \end{equations}
A direct computation then shows that the following laws of evolution hold
\begin{equation}\label{te2}
 \frac{du(t)}{dt}  = -u(t), \qquad \frac{d\theta(t)}{dt} = 2\left( 1- \theta(t)\right).
 \end{equation}
Of course, in the case in which $\varphi$ satisfies \fer{norm}, \fer{te2} imply conservation of both mean velocity and temperature.

It is well-known that the solution to the Fokker--Planck equation \fer{FP} converges exponentially in time to zero in relative entropy \cite{Tos97,Tos99}, which implies exponential convergence to equilibrium in $L^1$-norm. A slightly less known result is that exponential convergence to equilibrium for non regular initial data can be directly shown to hold also in weaker norms by resorting to the Fourier transform. 
Given a probability density $f(v)$, $v \in \R^d$, we define its \emph{Fourier transform} $\widehat{f}(\xi)$,  $\xi \in \R^d$ by
    \begin{equation*}
    \widehat{f}(\xi) = \int_{\R^d}e^{- i \, \xi\cdot v} f(v)\, dv.
    \end{equation*}
 Let us consider a family of metrics that has been introduced in the paper \cite{GTW} to study the trend to equilibrium of solutions to the space homogeneous Boltzmann equation for Maxwell molecules, and subsequently applied to a variety of problems related to kinetic models of Maxwell type.
For a more detailed description, we address the interested reader to the  lecture notes \cite{Carrillo:2007}.

Given $s >0$ and two random variables $X,Y$ with probability distributions $f$ (respectively $g$), their Fourier based distance $d_s(X,Y)$ is given by the quantity
\begin{equation*}
d_s(X,Y) = d_s(f,g) := \sup_{\xi \in \R^d} \frac{\left |\widehat{f}(\xi)-\widehat{g}(\xi)\right |}{|\xi|^s}.
\end{equation*}
The distance is finite, provided that $X$ and $Y$ have the same moments up to order $[s]$, where, if $s \notin \NN$,  $[s]$ denotes the entire part of $s$, or up to order $s-1$ if $s \in \NN$. Moreover $d_s$ is an ideal metric.  Its main properties are the following
\begin{enumerate}
\item Let $X_1,X_2,X_3$, with $X_3$ independent of the pair $X_1,X_2$ be  random variables with probability distributions $f_1$, $f_2$, $f_3$. Then
 \[ 
d_s(X_1 + X_3, X_2+X_3) = d_s(f_1* f_3, f_2*f_3) \leq d_s (f_1, f_2)= d_s(X_1, X_2) 
 \]
where the symbol $*$ denotes convolution;
\item Define, for a given nonnegative constant $a$, the dilation of a function $f= f(v)$, $v \in \R^d$ as
 \[
 f_a(v) = \frac 1{a^d} f\left ( \, \frac va \, \right ).
  \] 
  Then, given two random variables $X,Y$ with probability distributions $f$ and $g$, for any nonnegative constant $a$
 \[
  d_s(aX, aY) = d_s( f_{a}, g_{a}) \leq a^s \, d_s(f, g)= a^s \, d_s(X,Y). 
  \]
\end{enumerate}
Consider that, in view of the representation formula \fer{sol}, given the initial values $\varphi$ and $\tilde\varphi$, the two solutions $f(t)$ and $\tilde f(t)$ are the probability distributions of the random variables $Z_t$ and $\tilde Z_t$ expressed  by \fer{sol}. Hence if for some $s>0$ the distance $d_s(X, \tilde X) = d_s(\varphi, \tilde\varphi)$ is finite, then
 \begin{equations}\label{ok}
 d_s(Z_t, \tilde Z_t) = &\,\, d_s(e^{-t}X + (1-e^{-2t})^{1/2}W, e^{-t}\tilde X + (1-e^{-2t})^{1/2}W) \le \\
 &d_s(e^{-t}X, e^{-t}\tilde X) \le e^{-st}d_s( X, \tilde X).
 \end{equations}
The first inequality follows from property $1$. Then the dilation property $2$ is applied to conclude.
The same result \cite{cato2} can be easily obtained also resorting to the Fourier transform version of the Fokker--Planck equation \fer{FFP}, that reads
  \be\label{FFP}
  \frac{\partial \ff(\xi,t)}{\partial t} = -|\xi|^2 \ff(\xi,t) - \sum_{k=1}^d \xi_k \frac{\partial \ff(\xi,t)}{\partial \xi_k }.
  \ee

\section{The one-dimensional Fokker--Planck equation}\label{one-d}

For the rest of this Section, let us fix $d=1$.  To start with, let us consider that the Fokker--Planck equation \fer{FP} can be fruitfully written in weak form. It corresponds to say that, for any given smooth function $\phi(v)$, the Fokker--Planck operator modifies the solution $f(v,t)$ according to
 \be
  \label{kine-w}
\frac{d }{d\tau}\int_{\R} \phi(v) f(v,t)\,dv = \int_{\R}
 \bigl(\phi''(v) - v\phi'(v)\bigr) f(v,t)  \, dv.
 \ee
 By choosing $\phi(v) = e^{-i\xi v}$ we obtain the (one-dimensional) Fourier transform version \fer{FFP} of the Fokker--Planck equation \fer{FP}. The advantage of working with a weak version of the equation, is that we can allow the initial value $\varphi(v)$ to be a measure on $\R$. 

\subsection{A Rosenau-type approximation}

Rosenau \cite{Ro92}  proposed a regularized version of the Chapman-Enskog expansion of hydrodynamics. This regularized expansion resembles the usual Navier-Stokes viscosity terms at law wave-numbers, but unlike the latter, it has the advantage of being a bounded macroscopic approximation to the linearized collision operator.
The model originally considered by Rosenau is given by the scalar equation
 \be\label{Rose}
 f_t +  \Psi(g)_v = \left[ \frac {- \ve\xi^2}{1 + \ve^2 \xi^2}\ff(\xi) \right]^{\vee},
 \ee
where $\ff(\xi)$ denotes the Fourier transform of $f(v)$, while $ f(\xi)^\vee $ denotes the inverse Fourier transform.

The operator on the right side looks like the usual viscosity term $\ve f_{vv}$ at low wave-numbers $\xi$, while for higher wave numbers it is intended to model a bounded approximation of a linearized collision operator, thereby avoiding the artificial instabilities that occur when the Chapman-Enskog expansion for such an operator is truncated after a finite number of terms.

The right side of \fer{Rose} can be written as
 \be\label{Ros2}
\left[ \frac {- \ve \xi^2}{1 + \ve^2 \xi^2}\ff(\xi) \right]^{\vee} = \frac 1\ve \left[ \frac 1 {1 +\ve^2  \xi^2}\ff(\xi) -\ff(\xi) \right]^{\vee} = \frac 1{ \ve} \left[ M_{\ve}*f -f\right],
 \ee
where  $*$ denotes convolution and
 \be \label{Max}
 M_\gamma(v) = \frac 1{2\gamma}\,e^{-|v|/\gamma} = \frac 12\left(\frac 1{\gamma}\,e^{-v/\gamma}I_{\{v \ge 0\}}(v) + \frac 1{\gamma}\,e^{v/\gamma}I_{\{v < 0\}}(v) \right) 
 \ee
 is a non-negative symmetric function satisfying $\| M_\gamma\|_{L^1} = 1$. In \fer{Max} $I_A(v)$ denotes the characteristic function of the set $A$, i.e. $I_A(v) =1$ if $v \in A$, while $I_A(v) = 0$ otherwise. 

Hence Rosenau approximation consists in substituting  the linear diffusion equation
 \be\label{heat}
\frac{\partial g(v,t)}{\partial t} = \frac{\partial^2 g(v,t)}{\partial v^2}
 \ee
with the linear kinetic equation
 \be\label{kin}
\frac{\partial g(v,t)}{\partial t} = \frac {1}{\ve^2} \left[ M_{\ve}*g(v,t) -g(v,t)\right]
 \ee
A detailed study of the properties of the kinetic equation \fer{kin}, as well as its connections with the diffusion equation \fer{heat} was recently given in \cite{RT} (cf. also \cite{LT, PT13}). Also, a similar approximation was used in connection with the one-dimensional fractional diffusion equation \cite{FPTT}.
These studies showed that the approximation maintains most of the properties of the original diffusion equations. Moreover, in the case of the heat equation, the Rosenau-type approximation gives a new physical inside into the numerical approximation of \fer{heat}. In particular the moments at the first two order of the approximate solutions follow the same evolution of the original diffusion equation. 

However, it is important to notice that \fer{Ros2} is only a physically relevant way to write the right-hand side of \fer{Rose}. Indeed, by making use of standard properties of the convolution operation one has
 \be\label{Ros3}
\left[ \frac {- \ve \xi^2}{1 + \ve^2 \xi^2}\ff(\xi) \right]^{\vee} = -\ve \xi^2  \left[ \frac 1{1 + \ve^2 \xi^2}\ff(\xi) \right]^{\vee} = \ve \frac {\partial^2 M_{\ve}*f}{\partial v^2} =\ve M_\ve * \frac {\partial^2 f}{\partial v^2}.
 \ee
This shows that Rosenau approximation is obtained by smoothing out the right-hand side of the linear diffusion equation \fer{heat} by means of its convolution with $M_\ve$. Hence, the linear kinetic equation \fer{kin} can be alternatively written as
 \be\label{kin2}
\frac{\partial g(v,t)}{\partial t} =  M_\ve* \frac{\partial^2 g(v,t)}{\partial v^2}.
 \ee
It is tempting to use the same approximation for the one-dimensional Fokker--Planck equation \fer{FP}.  In this case one considers the equation
\be\label{FP-aa}
\frac{\partial f(v,t)}{\partial t} =  M_\ve*\left[ \frac{\partial^2 f(v,t)}{\partial v^2} +  \frac{\partial(v f(v,t))}{\partial v}\right] .
 \ee
In Fourier transform, equation \fer{FP-aa} reads
 \be\label{FFP1}
  \frac{\partial \ff(\xi,t)}{\partial t} = -\frac{\xi^2}{1 +\ve^2\xi^2} \ff(\xi,t) - \frac{\xi}{1 +\ve^2\xi^2} \frac{\partial \ff(\xi,t)}{\partial \xi }.
  \ee
It is clear that the diffusion term is given by a linear kinetic equation of type \fer{Ros2}. Then,  the Fourier transform of the drift term can be written as
 \[
 - \frac{\xi}{1 +\ve^2\xi^2} \frac{\partial \ff(\xi,t)}{\partial \xi } = \frac 1{2\ve} \left[ \frac 1{1-i\ve \xi} -  \frac 1{1+ i\ve \xi}\right]\frac{\partial(i \ff(\xi,t))}{\partial \xi }. 
 \]
Hence we obtain the identity 
 \be\label{drift}
M_\ve* \frac{\partial(v f(v,t))}{\partial v} = \frac 1\ve \tilde M_\ve * (vf(v,t)),
 \ee
where 
 \be \label{Max5}
 \tilde M_\gamma(v) = \frac 12\left(\frac 1{\gamma}\,e^{v/\gamma}I_{\{v < 0\}}(v) - \frac 1{\gamma}\,e^{-v/\gamma}I_{\{v \ge 0\}}(v) \right). 
 \ee
Note that, while $M_\gamma(v)$ is a symmetric probability density, $\tilde M_\gamma$ is antisymmetric, and it is obtained from $M_\gamma(v)$ by changing the sign on the domain $v \ge 0$. 

Finally, the approximated Fokker--Planck equation \fer{FP-aa} can be equivalently written as
 \be\label{Ros11}
 \frac{\partial f(v,t)}{\partial t} = \frac {1}{\ve^2} \left[ M_{\ve}*f(v,t) -f(v,t)\right] + \frac 1\ve\tilde M_\ve * (vf(v,t)). 
 \ee
One can easily verify, in view of the  expressions of the functions $M_\ve$ and $\tilde M_\ve$, that the mean velocity and temperature of the solution to \fer{Ros11} follow the same laws of evolution \fer{te2} of the original Fokker--Planck equation \fer{FP}.

\subsection{The discrete Fokker--Planck equation}\label{FFPPP}
Having in mind the previous discussion, given a small positive parameter $\ee\ll 1$ we consider the approximation to the operator $J_{FP}$ given by
\be\label{FPapp}
 J_{FP}^\ve(f)(v) = \frac 2{\ve^2}\left( P_{\ve} * f - f \right)(v)
  + \frac 1\ve \tilde P_\ve * (vf(v)),
  \ee
acting on probability densities $f(v)$ satisfying the normalization conditions \fer{norm}.  
In \fer{FPapp} 
 \begin{equation}\label{ker}
 P_\ve(v) = \frac 12\left( \delta(v+\ve) + \delta (v-\ve)\right), \quad
 \tilde P_\ve(v) = \frac 12\left( \delta(v+\ve) - \delta (v-\ve)\right),
 \end{equation}  
where $\delta(v)$ denotes as usual the Dirac delta function concentrating on $v=0$. Note that, in analogy with \fer{Ros11}, $P_\ve(v)$ is a symmetric probability measure, and $\tilde P_\ve$ is antisymmetric, and obtained from $P_\ve$ by changing its sign when $v >0$. 

Using the symmetry properties of $P_\ve$ and $\tilde P_\ve$, one shows that the  weak form of the evolution equation
 \be\label{FP-deb}
 \frac{\partial f_\ve}{\partial t} = J_{FP}^\ve(f_\ve) 
 \ee
 is written as
  \be
  \label{app1}
\frac{d }{dt}\int_{\R} \phi(v) f_\ve(v,t)\,dv = \int_{\R}\left[ \frac 2{\ve^2}\left( P_\ve*\phi - \phi\right)(v) - v \, \frac 1\ve \tilde P_\ve*\phi(v) \right] f_\ve (v) \, dv.
 \ee
 Alternatively, by choosing $\Phi(v) = e^{-i\xi v}$ one obtains the evolution equation for the Fourier transform $\ff_\ve(\xi,t)$
 \be\label{FFPapp}
  \frac{\partial \ff_\ve(\xi,t)}{\partial t} = -\frac 2{\ve^2}(1-\cos(\ve\xi))\ff_\ve(\xi,t) - \frac{\sin(\ve\xi)}\ve \frac{\partial \ff_\ve(\xi,t)}{\partial \xi }.
 \ee
The Fourier description  clearly shows that as $\ve \to 0$, the right-hand side of equation \fer{FFPapp} converges pointwise to the right-hand side of equation \fer{FFP}. 

Let us examine in details the properties of the solution to equation \fer{app1}. In view of definition \fer{ker}, it follows that for $0\le n \in \N$
 \begin{equations}\label{mo1}
 &\int_\R v^{2n} P_\ve(v)\, dv  = \ve^{2n},\qquad \int_\R v^{2n+1} P_\ve(v)\, dv  = 0, \\
 &\int_\R v^{2n} \tilde P_\ve(v)\, dv  = 0, \qquad \int_\R v^{2n+1} \tilde P_\ve(v)\, dv  = -\ve^{2n+1}.
 \end{equations}
 Hence, by choosing $\phi(v) = 1, v, v^2$ in \fer{app1} and using \fer{mo1} we conclude that equation \fer{app1} preserves the total mass, while the evolution equations for the mean velocity and temperature coincide with equations \fer{te1}. 

In particular, if the initial datum satisfies the normalization conditions \fer{norm}, also the solution to \fer{FP-deb} satisfies \fer{norm}.

A further fundamental property of the solution to equation \fer{app1} follows by considering as initial datum the law of a random variable $X_\ve$ that takes values only on a discrete number of points. To be more precise, for a given positive number $N \in \N$, $N \gg 1$, let us set $\ve = 1/N$. Then we define
 \be\label{ini1}
 \varphi_\ve(v) = \sum_{|j|\le 2N^2}\varphi_j \delta (v - j\ve), \quad \varphi_\ve(v) = 0  \,\,\, \rm{elsewhere}
 \ee
where the nonnegative constants $\varphi_j$ satisfy
 \be\label{nor2}
\sum_{|j|\le 2N^2}\varphi_j =1. 
 \ee 
Let  $\mathcal D$ be the space of functions of type \fer{ini1}, subject to condition \fer{nor2}. Let us consider a nonnegative measure $g(v)$ in $\mathcal D$. Owing to definition \fer{ker}, it is immediate to verify that $T(g)= \ve^2 J_{FP}^\ve(g)/2+ g$ belongs to $\mathcal D$. Indeed
 \begin{equations}\label{vero}
 T(g)(v) = &\,\,  P_{\ee} * g(v) + \frac\ve{2} \tilde P_\ve * (vg(v))=\\
 &\frac 12 \left[1 + \frac \ve{2} (v +\ve) \right]g(v+\ve) + \frac 12 \left[1 - \frac \ve{2} (v - \ve) \right]g(v-\ve).
 \end{equations}
Therefore $T(g)(v)$ is a linear combination of the values of $g$ in the points $g(v+\ve)$ and $g(v-\ve)$, and whenever $v$ belongs to the interval $(-2/\ve, 2/\ve)$ the coefficients of the linear combination  are nonnegative. Moreover, a direct inspection shows that
 \be\label{bound}
 T(g)(2/\ve + \ve) = T(g)(-2/\ve - \ve) = 0.
 \ee
 Last, condition \fer{nor2} remains verified. This property can be verified  by owing to the mass conservation property of $J_{FP}^\ve$, or by direct inspection. Indeed, using \fer{bound} we have
 \[
\sum_{|j|\le 2N^2}T(g)_j = \frac 12  \sum_{|j|\le 2N^2}\left( g_{j+1} + g_{j-1}\right) +
\frac {\ve^2}4 \sum_{|j|\le 2N^2}\left((j+1) g_{j+1} - (j-1)g_{j-1}\right) =
 \]
 \[
 \frac 12  \sum_{j=-2N^2+1}^{2N^2-1}g_j + \frac 12 \left( g_{2N^2} + g_{-2N^2}\right) +
 \frac {\ve^2}4 \left( 2N^2g_{2N^2} + 2N^2g_{-2N^2} \right) = \sum_{|j|\le 2N^2}g_j.
 \]
 Hence
  \[
T(g)(v) = \sum_{|j|\le 2N^2}T(g)_j \delta (v - j\ve), \quad T(g)(v) = 0  \,\,\, {\rm{elsewhere}, \quad \rm{and}}\,\,\,{\sum_{|j|\le 2N^2}T(g)_j =1}. 
 \] 
 Consequently, $T(g)$ is a linear mapping of $\mathcal D$ into $\mathcal D$. 

Now, given the nonnegative measure $\varphi_\ve \in \mathcal D$, consider that the initial value problem for the discrete Fokker--Planck equation \fer{FP-deb} can be rewritten as 
 \be\label{FPW}
\frac{\partial f_\ve(v,t)}{\partial t} = \frac 2{\ve^2}\left( T(f_\ve)(v) -  f_\ve(v)\right), \quad f_\ve(v, t=0) = \varphi_\ve(v). 
 \ee
The (unique) solution to \fer{FPW} can be explicitly expressed in the form of a Wild sum  \cite{PT13,Wild}
 \be\label{solW}
 f_\ve(v,t) = e^{-2t/\ve^2}\sum_{i \ge 0} \frac 1{i!}\left(\frac{2t}{\ve^2} \right)^i f_\ve^{(i+1)}(v),
 \ee
 where $f_\ve^{(0)}(v) = \varphi_\ve(v)$ is the initial value, and the nonnegative coefficients $f_\ve^{(i)}$, $i \ge 1$, are recursively defined by $ f_\ve^{(i)}(v) = T(f_\ve^{(i-1)})(v)$. Since at any time $t \ge 0$ the solution \fer{solW} is a convex combination of the time-independent nonnegative coefficients $f_\ve^{(i)}$, the solution to the initial value problem \fer{FPW} with $\varphi_\ve \in \mathcal D$ belongs to $\mathcal D$ at any subsequent time $t \ge0$. In addition, the solution is nonnegative for all $t \ge 0$. 

Last, let $\Psi(r)$, $r \ge 0$ be a convex function, such that $\Psi(0) = 0$. Then,  given a nonnegative measure $g(v)$ in $\mathcal D$, by \fer{vero}, whenever $v$ belongs to the interval $(-2/\ve, 2/\ve)$,  we obtain the inequality
 \be\label{cov2}
\Psi(T(g)(v)) \le \frac 12 \left[1 + \frac \ve{2} (v +\ve) \right]\Psi\left(g(v+\ve)\right) + \frac 12 \left[1 - \frac \ve{2} (v - \ve) \right]\Psi\left(g(v-\ve)\right),
 \ee
while
 \[
 \Psi\left( T(g)(2/\ve + \ve)\right) = \Psi\left( T(g)(-2/\ve - \ve)\right) = 0.
 \]
Therefore, proceeding as before, we obtain
\be\label{ent3}
\sum_{|j|\le 2N^2}\Psi\left(g_j \right) \le  \sum_{|j|\le 2N^2}\Psi\left( T(g)_j\right).
\ee
Since at any time $t \ge 0$ the solution \fer{solW} is a convex combination of the time-independent coefficients $f_\ve^{(i)}$, the convexity of $\Psi$, coupled with inequality \fer{ent3} shows that, for any given nonnegative measure $\varphi_\ve \in \mathcal D$
 \be\label{decre}
\Psi(f)(t) = \sum_{|j|\le 2N^2}\Psi\left(f_{\ve,j}(t) \right) \le  \sum_{|j|\le 2N^2}\Psi\left( \varphi_j\right) = \Psi(\varphi_\ve).
 \ee
Hence, the discrete  functional $\Psi(f_\ve)(t)$ is monotonically decreasing in time. In particular, we can consider the classical Shannon entropy on the probability measure $g \in \mathcal D$, defined by
 \be\label{sha}
 H(g) = - \sum_{|j|\le 2N^2}g_j \log g_j. 
 \ee
Then, in analogy with the original Fokker--Planck equation, Shannon entropy \fer{sha} is shown to be increasing in time along the solution to the discrete Fokker--Planck equation \fer{FP-deb}.
  
A further interesting property of the approximation \fer{FP-deb} is the fact that the equation
 $J_{FP}^\ve(f) =0$ has a unique explicit solution in $\mathcal D$, which is nothing but the stationary solution of the  approximation \fer{app1}, provided that the initial measure $\varphi_\ve \in \mathcal D$.

\subsection{The stationary solution}\label{stazz}

In one dimension, the Fourier transform of the stationary solution solves the equation
 \be\label{statio}
\frac 2{\ve}(1-\cos(\ve\xi))\ff_\infty(\xi) + \sin(\ve\xi) \frac{\partial \ff_\infty(\xi)}{\partial \xi } = 0.
 \ee
Hence, if $\ff_\infty(\xi) \neq 0$
 \[
 \frac 1{\ff_\infty(\xi)}\frac{\partial \ff_\infty(\xi)}{\partial \xi } = - \frac 2{\ve}\,\frac{1-\cos(\ve\xi)}{\sin(\ve\xi)} = - \frac 2{\ve^2}\,\frac{\ve \sin(\ve\xi)}{1+\cos(\ve\xi)} = \frac 2{\ve^2}\, \frac d{d\xi} \log (1+ \cos(\ve\xi)).
 \] 
Therefore,  integrating both sides from $0$ to $\xi$, and assuming that the solution  has mass equal to one, so that $\ff_\infty(0) =1$ one obtains
 \be\label{sta5}
\ff_\infty(\xi) = \left( \frac{1+ \cos(\ve\xi)} 2\right)^{2/\ve^2}.
 \ee
By choosing  $ \ve = 1/N$, we argue that the steady state is the convolution product of $2N^2$ identical functions, each of them with Fourier transform
 \be\label{tra2}
 \widehat\psi(\xi) = \frac{1+ \cos{\frac \xi N}}2. 
 \ee
Let $X$ be a discrete random variable taking values $\pm 1$ and $0$ with probabilities
 \be\label{var6}
 P(X = \pm 1) = \frac 14, \qquad P(X=0) = \frac 12.
 \ee
Then, 
 \[
 h(x) = \frac 14 \left( \delta(x+1) + \delta(x-1) + 2 \delta (x) \right),
 \]
which implies
 \[
 \widehat h(\xi) = \frac{1+ \cos{ \xi}}2.
 \]
Hence, the function \fer{tra2} is the characteristic function of the random variable $X/N$. Consequently, if $ \ve = 1/N$, the function \fer{sta5} is the characteristic function of the random variable
 \be\label{NN}
 S_N = \frac 1 N \sum_{j=1}^{2N^2} X_j,
 \ee
where the $X_j$ are independent and identically distributed copies of $X$, defined as in \fer{var6}. By construction, the random variable $S_N$ takes values only on the discrete set of points $\ve j$, where $|j| \le 2N^2$. Consequently, $S_N$ has a probability distribution that belongs to $\mathcal D$. Considering now that $X$ has zero mean, and variance $1/2$, it follows easily by the central limit theorem that the law of $S_N$ is an approximation of the Gaussian distribution, that is an approximation of the stationary solution of the original Fokker--Planck equation. In addition, the law of $S_N$ satisfies the normalization conditions \fer{norm}.

\subsection{Large-time behavior}\label{larg}
The results of Section \ref{FFPPP} showed that the solution to the discrete Fokker--Planck equation \fer{FP-deb} in $\mathcal D$ maintains most of the properties of the original Fokker--Planck equation, like preservation of positivity and mass,  same evolution of moments up to order two, and entropy monotonicity. A further property of the original Fokker--Planck equation is the exponential convergence of its solution towards the Maxwellian equilibrium. As discussed in Section \ref{model}, exponential convergence to equilibrium can be shown also in Fourier metric.  
In what follows, we will investigate about the large-time behavior of the solution to \fer{FP-deb} and its (eventual) convergence to equilibrium in terms of the metric $d_s$.

Let $\ve\xi \in [m\pi, (m+1)\pi)$, where $|m| \in \N$. Then, since the function $\sin\ve\xi$ does not change sign in this interval, the differential equation
 \be\label{chara}
 \frac{ d\xi(t)}{dt} = \frac 1\ve \sin(\ve\xi(t)), \qquad \xi(t=0)=\xi,
 \ee
can be solved uniquely as soon as $\ve\xi(t)\in [m\pi, (m+1)\pi)$ to give the relationship
 \be\label{rel3}
 \tan \frac{\ve\xi(t)}2 = e^t\, \tan \frac{\ve\xi}2.
 \ee
Then, since both $\ve\xi(t)/2$ and $\ve\xi/2$ belong to the interval $[m\frac \pi 2, (m+1)\frac \pi 2)$ identity \fer{rel3} can be used to relate in a unique way $\xi$ to $\xi(t)$ (or vice-versa)
 \be\label{rel5}
  \xi(t) = \frac 2\ve \arctan\left[ e^t\, \tan \frac{\ve\xi}2 \right].
 \ee
Now, considering that 
 \[
 1 - \cos(\ve\xi) = \frac{2 \left(\tan (\ve\xi/2)\right)^2}{1 + \left(\tan (\ve\xi/2)\right)^2},
 \]
by using \fer{rel3}, the one-dimensional equation \fer{FFP} can be integrated along characteristics on each interval $[m\pi, (m+1)\pi)$. Indeed, using \fer{rel3} on this interval equation \fer{FFP} takes the form
 \begin{equations}\label{facile}
 &\frac d{dt} \ff_\ve(\xi(t),t) = -\frac 2{\ve^2} \frac{2 \left(\tan (\ve\xi(t)/2)\right)^2}{1 + \left(\tan (\ve\xi(t)/2)\right)^2} \ff_\ve(\xi(t),t)= \\
 &-\frac 2{\ve^2} \frac{2\,e^{2t} \left(\tan (\ve\xi/2)\right)^2}{1 +e^{2t} \left(\tan (\ve\xi/2)\right)^2} \ff_\ve(\xi(t),t) = -\frac 2{\ve^2}\frac d{dt} \log \left[1 +e^{2t} \left(\tan (\ve\xi/2)\right)^2\right] \ff_\ve(\xi(t),t).
 \end{equations}
Thus, integration over time gives the solution
 \be\label{sol5}
 \ff_\ve(\xi(t),t) = \widehat\varphi_\ve(\xi) \left[\frac{1 + \left(\tan (\ve\xi/2)\right)^2}{1 + \left(\tan (\ve\xi(t)/2)\right)^2}\right]^{2/\ve^2}.
 \ee
Then using again \fer{rel3} we obtain for $\xi \in [m\pi, (m+1)\pi)$ 
 \be\label{solu}
\ff_\ve(\xi,t) = \widehat\varphi_\ve\left[ \frac 2\ve \arctan\left(e^{-t}\tan (\ve\xi/2)\right)\right]\left[\frac{1 + e^{-2t}\left(\tan (\ve\xi/2)\right)^2}{1 + \left(\tan (\ve\xi/2)\right)^2}\right]^{2/\ve^2}
 \ee
We remark that,  if the initial value $\varphi_\ve(v)\in \mathcal D$, its Fourier transform is given by
 \be\label{ini-f}
 \widehat\varphi_\ve(\xi) = \sum_{|j|\le N^2}\varphi_j e^{-i\ve \xi j}.
 \ee
where the nonnegative constants $\varphi_j$ satisfy condition \fer{nor2}. Then, for any value $\bar\xi = 2m\pi/\ve$, with $m \in \N$, $ \widehat\varphi_\ve(\bar\xi) = 1$. By letting $t \to \infty$ in \fer{solu}, for any $\xi \in \R$ such that $\ve\xi \in [m\pi, (m+1)\pi)$
 \[
\widehat\varphi_\ve\left[ \frac 2\ve \arctan\left(e^{-t}\tan (\ve\xi/2)\right)\right] \to  \widehat\varphi_\ve\left(\frac 2{\ve}\, m\pi\right) = 1,
 \]
while
 \be\label{sta1}
 \ff_\ve(\xi,t) \to \left(\frac 1{1+ \left(\tan(\ve\xi/2)\right)^2}\right)^{2/\ve^2}=\left(\frac{1 + \cos(\ve\xi)} 2 \right)^{2/\ve^2} = \ff_\infty(\xi).
 \ee
Since $m$ is arbitrary, pointwise convergence to the stationary solution \fer{sta5} follows for all $\xi \in \R$. 

A stronger result about convergence to the steady state follow by restricting the allowed set of values of $\xi$. 

Indeed, since the function \fer{sta5} is a stationary solution to the Fokker--Planck equation, for every $t \ge0$ it satisfies the identity
\be\label{solud}
\ff_\infty(\xi) = \ff_\infty\left[ \frac 2\ve \arctan\left(e^{-t}\tan (\ve\xi/2)\right)\right]\left[\frac{1 + e^{-2t}\left(\tan (\ve\xi/2)\right)^2}{1 + \left(\tan (\ve\xi/2)\right)^2}\right]^{2/\ve^2}.
 \ee
Therefore, considering that for $t \ge 0$
 \[
 \left[\frac{1 + e^{-2t}\left(\tan (\ve\xi/2)\right)^2}{1 + \left(\tan (\ve\xi/2)\right)^2}\right]^{2/\ve^2} \le 1,
 \]
if the initial value $\varphi_\ve(v)$ has zero mean, we have the inequality
 \[
 \frac{|\ff_\e(\xi,t) -\ff_\infty(\xi) |}{|\xi|^2} \le \frac{\left|\widehat\varphi_\ve -\ff_\infty \right|\left( \frac 2\ve \arctan\left(e^{-t}\tan (\ve\xi/2)\right) \right)}{|\xi|^2},
 \]
and this inequality, on the set $|\ve\xi| \le \pi/2$, clearly implies 
 \be\label{dec3}
 \sup_{|\ve\xi| \le \pi/2}\frac{|\ff_\e(\xi,t) -\ff_\infty(\xi) |}{|\xi|^2}  \le d_2(\varphi_\ve,f_\infty) \sup_{|\ve\xi| \le \pi/2} \left|\frac  {\arctan\left(e^{-t}\tan (\ve\xi/2)\right)}{\frac{\ve\xi}2}\right|^2.
 \ee
Hence, by setting $\eta = \tan(\ve\xi/2)$ we obtain
 \[
 \sup_{|\ve\xi| \le \pi/2} \left|\frac  {\arctan\left(e^{-t}\tan (\ve\xi/2)\right)}{\frac{\ve\xi}2}\right|^2 = \sup_{|\eta| \le 1} \left|\frac  {\arctan\left(e^{-t}\eta \right)}{\arctan \eta}\right|^2 \le 4 e^{-2t}.
 \]
This implies exponential convergence of the solution to the stationary state (on the set $|\ve\xi| \le \pi/2$) at the same rate of the Fokker--Planck equation \fer{FP}. On the other hand, on the set $|\ve\xi| > \pi/2$, since both $|\ff_\e|(\xi,t)\le1 $ and $|\ff_\infty|(\xi) \le 1$, we have the bound
 \be\label{bb7}
  \sup_{|\ve\xi| > \pi/2}\frac{|\ff_\e(\xi,t) -\ff_\infty(\xi) |}{|\xi|^2}  \le  \sup_{|\ve\xi| > \pi/2}\frac{2}{|\xi|^2} = \frac 8{\pi^2}\, \ve^2.
 \ee

\subsection{Stability of the approximation}
Let us consider a Taylor expansion of the trigonometric functions on the right-hand side of equation \fer{FFPapp} up to the second order. We obtain
 \be\label{tay}
  \frac{\partial \ff_\ve(\xi,t)}{\partial t} = -|\xi|^2 \ff_\ve(\xi,t) - \xi\frac{\partial \ff_\ve(\xi,t)}{\partial \xi} + R_\ve(\ff_\ve)(\xi,t), 
 \ee
where the remainder term has the expression 
 \be\label{rest}
  R_\ve(\ff_\ve)(\xi,t) = \xi^3 \left( \frac{2\ve}{3!}\sin \ve \bar\xi \, \ff_\ve(\xi,t) - \frac{\ve^2}{3!} \cos \ve\bar\xi\, \frac{\partial \ff_\ve(\xi,t)}{\partial \xi}\right),
 \ee
 and $\bar\xi$ belongs to the interval $(0, \xi)$. 
Let the initial value $\varphi$ of the Fokker--Planck equation \fer{FP} possess finite moments up to the order three, and let us consider an approximation $\varphi_\ve \in \mathcal D$  with the same moments of $\varphi$ up to the second order. Then, since the mean velocity and the temperature of the Fokker--Planck equation and of its approximation follow the  laws \fer{te2}, it is immediate to conclude that 
 \be\label{bbb}
 \frac{|R_\ve(\ff_\ve)(\xi,t)|}{|\xi|^3} \le  \frac{2\ve}{3!}\sup_\xi |\ff_\ve(\xi,t)| + \frac{\ve^2}{3!} \sup_\xi \left| \frac{\partial \ff_\ve(\xi,t)}{\partial \xi}\right| \le \frac{2\ve}{3!} + \frac{\ve^2}{3!}|u_0| = \ve C(\ve,u_0),
 \ee
where $u_0$ is the initial mean velocity defined by \fer{velo}. In addition, if the Fourier metric $d_3(f,f_\ve)(t)$ is initially bounded, it remains bounded at any subsequent time.

Let $h(\xi,t)$ be defined as
 \[
 h(\xi,t) = \frac{\ff(\xi,t) - \ff_\ve(\xi,t)}{|\xi|^3}. 
 \]
Since for $\xi \neq 0$ we have the identity
 \[
 \frac \xi{|\xi|^3} \frac{\partial(\ff(\xi,t)- \ff_\ve(\xi,t))}{\partial \xi} = 3\, h(\xi,t) + \xi  \frac{\partial h(\xi,t)}{\partial \xi}
 \]
by considering the difference between the Fourier transforms of one-dimensional Fokker--Planck equation \fer{FFP} and \fer{tay},  we conclude that $h$ satisfies 
 \be\label{h2}
  \frac{\partial h(\xi,t)}{\partial t} + \xi  \frac{\partial h(\xi,t)}{\partial \xi} = -(|\xi|^2 +3)h(\xi,t) + \frac{R_\ve(\ff_\ve)(\xi,t)}{|\xi|^3} .
 \ee 
Integrating along characteristics,  \fer{h2} is equivalent to
 \[
  \frac{d h(\xi e^t,t)}{d t}  = -(|\xi|^2 e^{2t} +3)h(\xi e^t,t) + \frac{R_\ve(\ff_\ve)(\xi e^t,t)}{|\xi|^3e^{3t}}.
 \]
Hence, by using \fer{bbb} we obtain
 \[
   \frac{d}{dt} |h(\xi e^t,t)|  \le -3 |h(\xi e^t,t)| + \ve C(\ve,u_0),
 \]
that implies 
 \[
 |h(\xi e^t,t)| \le |h_0(\xi)|e^{-3t} + \ve C(\ve,u_0)\left( 1- e^{-3t} \right).
 \]
In conclusion
 \be\label{fin}
 d_3(f(t), f_\ve(t)) = \sup_{\xi \in \R} |h(\xi,t)| \le  d_3(\varphi, \varphi_\ve) e^{-3t} + \ve C(\ve,u_0)\left( 1- e^{-3t} \right).
 \ee
Hence, by choosing initial data for the discrete Fokker--Planck equation such that $d_3(\varphi, \varphi_\ve) \le C\ve$, we obtain that in the Fourier distance $d_3$ the uniform in time estimate
 \be\label{unif}
d_3(f(t), f_\ve(t)) \le \ve \max\left\{ C, C(\ve,u_0)\right\}.
 \ee

\section{The general case}
\setcounter{equation}{0}

The discretization of the one-dimensional Fokker--Planck equation introduced in Section \ref{one-d} can be easily extended to any dimension $d>1$. To this aim, it is enough to outline that the Fourier transformed equation \fer{FFP} can be rewritten as
  \be\label{FFd}
  \frac{\partial \ff(\xi,t)}{\partial t} = -\sum_{k=1}^d\left( |\xi_k|^2 \ff(\xi,t) + \xi_k \frac{\partial \ff(\xi,t)}{\partial \xi_k }\right).
  \ee
The Rosenau approximation in this case is given by
 \be\label{Rosd}
  \frac{\partial \ff(\xi,t)}{\partial t} = -\sum_{k=1}^d\left(\frac{\xi_k^2}{1 + \ve^2\xi_k^2} \ff(\xi,t) + \frac{\xi_k}{1 + \ve^2\xi_k^2}  \frac{\partial \ff(\xi,t)}{\partial \xi_k }\right),
  \ee
that in the physical space reads
 \be\label{Ros11d}
 \frac{\partial f(v,t)}{\partial t} = \frac {1}{\ve^2}\sum_{k=1}^d\left( \left[ M_{\ve}(v_k)*f(v,t) -f(v,t)\right] + \frac 1\ve\tilde M_\ve(v_k) * (vf(v,t))\right). 
 \ee
As in Section \ref{one-d}, one can easily verify that the mean velocity and temperature of the solution to \fer{Ros11} follow the same laws of evolution of the original Fokker--Planck equation \fer{FP}.
Also, using the functions \fer{ker} one obtains the evolution equation for the Fourier transform $\ff_\ve(\xi,t)$
 \be\label{FFdapp}
  \frac{\partial \ff_\ve(\xi,t)}{\partial t} = -\sum_{k=1}^d\left(\frac 2{\ve^2}(1-\cos(\ve\xi_k))\ff_\ve(\xi,t) + \frac{\sin(\ve\xi_k)}\ve \frac{\partial \ff_\ve(\xi,t)}{\partial \xi_k }\right).
 \ee
In the limit $\ve \to 0$, the right-hand side of equation \fer{FFPapp} converges pointwise to the right-hand side of equation \fer{FFd}. 

As in the one-dimensional case, let us consider as initial datum   the law  $\varphi_\ve(v)= \varphi(v_1, v_2, \dots, v_d)$ of a random vector $(X^{(1)}_\ve, X^{(2)}_\ve, \dots, X^{(d)}_\ve)$ which takes values only on a discrete number of points of $\R^d$. For a given positive number $N \in \N$, $N \gg 1$, let us set $\ve = 1/N$, the law reads
 \be\label{inid}
 \varphi_\ve(v) = \sum_{i=1}^d \sum_{|j_i|\le 2N^2}\varphi_{j_1,\dots,j_d} \prod_{i=1}^d\delta (v_i - j_i\ve), \quad \varphi_\ve(v) = 0  \,\,\, \rm{elsewhere}
 \ee
where the nonnegative constants $\varphi_{j_1,\dots,j_d}$ satisfy
 \be\label{nor2d}
\sum_{i=1}^d\sum_{|j_i|\le 2N^2}\varphi_{j_1,\dots,j_d} =1. 
 \ee 
Let  $\mathcal D^d$ be the space of functions of type \fer{ini1}, subject to condition \fer{nor2}. Proceeding as in Section \ref{one-d}, we can prove that, starting from a nonnegative initial value in $\mathcal D^d$, the solution to equation \fer{Ros11d} belongs to $\mathcal D^d$ for each time $t \ge 0$. Moreover, the solution is nonnegative. Also, mass is preserved, and the laws of evolution of the mean velocity and temperature follow the same laws of evolution \fer{te2} of the continuous equation. Last, Shannon entropy is monotonically increasing. Concerning the equilibrium distribution, we assume that
 \be\label{sta5d}
\ff_\infty(\xi) = \prod_{j=1}^d \ff_\infty(\xi_j)= \prod_{j=1}^d \left( \frac{1+ \cos(\ve\xi_j)} 2\right)^{2/\ve^2}.
 \ee
Then we have
\begin{equations}\label{sta-d}
& \sum_{k=1}^d\left(\frac 2{\ve^2}(1-\cos(\ve\xi_k)) \ff_\infty(\xi) +  \frac{\sin(\ve\xi_k)}\ve \frac{\partial \ff_\infty(\xi)}{\partial \xi_k }\right)= \\
&\sum_{k=1}^d\prod_{j\not=k}\ff_\infty(\xi_j)\left(\frac 2{\ve^2}(1-\cos(\ve\xi_k)) \ff_\infty(\xi_k) +  \frac{\sin(\ve\xi_k)}\ve \frac{\partial \ff_\infty(\xi_k)}{\partial \xi_k }\right) = 0.
  \end{equations}
Consequently, $\ff_\infty(\xi)$ is a stationary solution for the $d$-dimensional Fokker--Planck equation. Note that 
 \[
\widehat\Psi(\xi) = \prod_{j=1}^d \left( \frac{1+ \cos(\ve\xi_j)} 2\right)
 \]
is the Fourier transform of the joint distribution function of the random vector $X=(X_1, \dots,X_d) $ where the random variables $X_k$, for $k=1,\dots, d$ are independent each other and distributed according to \fer{var6}.
Consequently, if $ \ve = 1/N$, the function \fer{sta5d} is the characteristic function of the random vector
 \be\label{NNd}
 S_N = \frac 1 N \sum_{j=1}^{2N^2} X_j,
 \ee
where the $X_j$ are independent and identically distributed copies of the random vector $X$ Considering now that $X$ has zero mean, and variance $d/2$, it follows by the central limit theorem that the law of $S_N$ is an approximation of the $d$-dimensional Gaussian distribution, that is an approximation of the stationary solution of the original Fokker--Planck equation in $\R^d$.

Unlikely, the analysis of Section \ref{larg} is no more valid in dimension $d>1$, and the study of the large-time behavior of the solution to the discretized Fokker-Planck equation requires further efforts. We leave this problem open for future research.

\section{Higher-order diffusions}
\setcounter{equation}{0}

Let us consider the fourth-order (one-dimensional) linear diffusion \cite{BG}
 \be\label{d4}
\frac{\partial g(v,t)}{\partial t} = - \frac{\partial^4 g(v,t)}{\partial v^4}.
 \ee
The Fourier transform version of equation \fer{d4} reads
 \be\label{f4}
\frac{\partial \widehat g(\xi,t)}{\partial t} =  -\xi^4 \widehat g(\xi,t).
 \ee
Following the idea of Rosenau \cite{Ro92}, for a given $\ve\ll 1$ we consider the approximation given by
  \be\label{f41}
\frac{\partial \widehat g(\xi,t)}{\partial t} =  -\left(\frac{\xi^2}{1+ \ve^2\xi^2}\right)^2 \widehat g(\xi,t)=  -\frac 1{\ve^4}\left(\frac{\ve^2\xi^2}{1+ \ve^2\xi^2}\right)^2 \widehat g(\xi,t).
 \ee
This approximation is consistent with the analogous one introduced for the Fokker--Planck equation. Note that, since 
 \[
 \left(\frac{1}{1+ \ve^2\xi^2}\right)^2 = \widehat{M_\ve*M_\ve}(\xi),
 \]
 with $M_\ve$ defined as in \fer{Max}, the approximation \fer{f41} corresponds to modify equation \fer{d4} by taking the convolution of the right-hand side with $M_\ve*M_\ve$. Hence the approximation in the physical space reads
 \be\label{ad4}
\frac{\partial g(v,t)}{\partial t} = -M_\ve*M_\ve*\frac{\partial^4 g(v,t)}{\partial v^4}=  \frac{\partial^4M_\ve*M_\ve*g(v,t)}{\partial v^4}.
 \ee

Since
 \[
\left(\frac{\ve^2\xi^2}{1+ \ve^2\xi^2}\right)^2 = \left(1- \frac{1}{1+ \ve^2\xi^2} \right)^2= 1 + \left(\frac{1}{1+ \ve^2\xi^2}\right)^2 -2 \frac{1}{1+ \ve^2\xi^2}. 
 \]
equation \fer{ad4} can be rewritten as a linear kinetic equation of the form
 \be\label{kin4}
\frac{\partial g(v,t)}{\partial t} = \frac {1}{\ve^4} \left[ G(M_\ve)*g(v,t) -g(v,t)\right],
 \ee
where the function
 \be\label{def-g}
  G(M_\ve)(v) = 2 M_\ve(v) - M_\ve*M_\ve(v).
 \ee
has integral equal to one, but, at difference with the case of the heat equation, is no more a probability density function. Indeed, it becomes negative on part of the real line. This fact is in agreement with the theory of higher-order diffusions, and it is connected with the absence of a maximum principle for the solution to these equations.

It is remarkable that, in view of the expression \fer{ad4}, the moments up to order four of the solution to the approximated equation  follow the same evolution of the original equation. 

The same property is maintained by considering in the definition of $G_\ve$ a probability density different from $M_\ve$. In particular, we can resort to $P_\ve$, as defined in \fer{ker}. Consequently, we consider the approximation
  \be\label{kin41}
\frac{\partial g(v,t)}{\partial t} = \frac {4}{\ve^4} \left[ G(P_\ve)*g(v,t) -g(v,t)\right],
 \ee
where the constant $4$ in front of the interaction operator is chosen to preserve the evolution of the fourth-order moment. By resorting to the definition of $P_\ve(v)$, it is immediate to conclude that 
\[
P_\ve*P_\ve(v) = \frac 14\left( \delta(v+2\ve) + 2\delta(v) + \delta(v-2\ve)\right).
\]
Hence, for a given (smooth) function $h=h(v)$
\be\label{dc4}
\frac {1}{\ve^4} \left[ G(P_\ve)*h(v) -h(v)\right] = \frac{ \delta(v+2\ve) -4 \delta(v+\ve)+ 6\delta(v) - 4\delta(v-\ve) + \delta(v-2\ve)}{\ve^4}. 
\ee
Expression \fer{dc4} coincides with one of the central differences approximation of the fourth order derivative of a function (cf. \cite{Lele} and the references therein). The previous reasoning allows to conclude that, at difference with other approximations, \fer{dc4} is well-adapted, in view of its properties about preservation of moments evolution, to approximate the diffusion equation \fer{d4}.

The approximation of the (linear) diffusion equation of order $2n$, with $n>2$ follows along the same lines. Indeed, given the equation
 \be\label{d2n}
\frac{\partial g(v,t)}{\partial t} = (-1)^{n-1} \frac{\partial^{2n} g(v,t)}{\partial v^{2n}},
 \ee
in Fourier variables, equation \fer{d2n} takes the form
\be\label{f2n}
\frac{\partial \widehat g(\xi,t)}{\partial t} =  -\xi^{2n} \widehat g(\xi,t).
 \ee
 Given $\ve\ll 1$, we consider the approximation to \fer{f2n} given by
  \be\label{f2n1}
\frac{\partial \widehat g(\xi,t)}{\partial t} =  -\left(\frac{\xi^2}{1+ \ve^2\xi^2}\right)^n \widehat g(\xi,t)=  -\frac 1{\ve^{2n}}\left(\frac{\ve^2\xi^2}{1+ \ve^2\xi^2}\right)^n \widehat g(\xi,t).
 \ee
In this case
\be\label{aba}
\left(\frac{\ve^2\xi^2}{1+ \ve^2\xi^2}\right)^n = \left(1- \frac{1}{1+ \ve^2\xi^2} \right)^n= 1 + \sum_{k=1}^n(-1)^k{n \choose k} \left(\frac{1}{1+ \ve^2\xi^2}\right)^k.
\ee
We can easily find the expression of the function on the right-hand side in the physical space by considering that, for $k \in \N$
 \[
\left(\frac{1}{1+ \ve^2\xi^2}\right)^k = \underbrace{\widehat{M_\ve*\cdots*M_\ve}}_{k}(\xi)
\]
Thus, the approximation to the linear diffusion equation of order $2n$, with $n \ge 2$,  can be rewritten as a linear kinetic equation of the form
 \be\label{kin2n}
\frac{\partial g(v,t)}{\partial t} = \frac {1}{\ve^{2n}} \left[ G_n(M_\ve)*g(v,t) -g(v,t)\right],
 \ee
where the function $G_n(M_\ve)$ is defined by
 \be\label{def-gn}
  G_n(M_\ve)(v) = \sum_{k=1}^n(-1)^{k+1}{n \choose k}\underbrace{{M_\ve*\cdots*M_\ve}}_{k}.
 \ee
As before, for $n \ge 2$
 \[
 \int_\R  G_n(M_\ve)(v)\, dv = 1.
  \]
 By construction, the solution to the kinetic equation \fer{kin2n} is such that, for all $k \le 2n-1$
 \[
  \int_\R  v^k\left[ G_n(M_\ve)*g(v,t) -g(v,t)\right] \, dv = 0,
   \]
 while
  \[
 \int_\R  v^{2n} \frac {1}{\ve^{2n}}\left[ G_n(M_\ve)*g(v,t) -g(v,t)\right] \, dv = (2n)! \int_\R g(v,t)\, dv.
  \]
Likewise, given a (smooth) function $h(v)$, the expression
 \be\label{cent2n}
C_{2n}(h)(v) =  \frac {1}{\ve^{2n}} \left[ G_n(P_\ve)*h(v) -h(v)\right] ,
 \ee
where as usual $P_\ve(v)$ is given by  \fer{ker}, gives an explicitly computable central difference approximation of order $2n$ of $h$ with a number of good properties with respect to the higher-order diffusion equation \fer{f2n}. Indeed, the solution to the approximated diffusion equation
 \[
 \frac{\partial g(v,t)}{\partial t} = (2n)! \frac {1}{\ve^{2n}} \left[ G_n(P_\ve)*g(v,t) -g(v,t)\right]  
  \]
is such that, in agreement with the solution to the original diffusion equation, all moments up to the order $2n-1$ remain constant in time, while the moment of order $2n$ is linearly increasing with the same rate. 

\section{Conclusions}

We introduced and discussed a discretized version of the Fokker--Planck equation which maintains most of the properties of the continuous version. The basic idea was to use a suitable modification of the Fourier transform of the equation, similar to the one considered by Rosenau \cite{Ro92} for the linear heat equation. This approach has been subsequently applied to higher-order linear diffusion operators, to obtain an easy-to-handle way to recover explicitly a central difference approximation to derivatives of any even order. A main problem, however, remains open. It is not clear wether the solution to the discrete Fokker--Planck equation converges towards the stationary discrete solution or not. A partial result in this direction has been derived in Section \ref{larg}. Also, it would be interesting to know if, in the set $\mathcal D$, and for random variables with a law that satisfies the normalization condition \fer{norm}, the Shannon entropy \fer{sha} attains the maximum value in correspondence to the law of the stationary distribution $S_N$ defined in \fer{NN}. 
\vskip 1cm
\noindent{\bf{Acknowledgments:}} This work has been written within the
activities of the National Group of Mathematical Physics (GNFM) of INdAM (National Institute of
High Mathematics), and partially supported by the   MIUR-PRIN Grant 2015PA5MP7 ``Calculus of Variations''.

\end{document}